\documentclass{article}

\usepackage{amsmath, amsfonts,amssymb, amsbsy}

\newtheorem{theo}{\bf Theorem}[section]
\newtheorem{cor}{\bf Corollary}[section]
\newtheorem{defi}{\bf Definition}[section]
\newtheorem{nota}{\bf Notation}[section]
\newtheorem{lem}{\bf Lemma}[section]

\begin{document}
\begin{center}
 {\Large   Perfect Graeco-Latin balanced incomplete
block designs \\ and related designs}
\end{center}

\vskip10pt
\begin{center}
  {\large Sunanda Bagchi\\
   Theoretical Statistics and Mathematics Unit  \\
Indian Statistical Institute\\ Bangalore 560059, India.}
\end{center}
\vskip10pt

{\bf Abstract :} Main effect plans orthogonal through the block
factor (POTB) have been defined and a few series of them have been
constructed in Bagchi (2010). These plans are very closely related
to the `mutually orthogonal balanced nested row-column designs' of
Morgan and Uddin (1996) and many other  combinatorial designs  in
the literature with different names like `BIBDs for two sets of
treatment', `Graeco-Latin designs' and `PERGOLAs'. In fact all of
them may be viewed as POTBs  satisfying one or more additional
conditions, making them `optimal'.  However, the PERGOLAs are
defined to satisfy an additional property, without which also it is
optimal. Interestingly, this additional property is satisfied by all
the hitherto known examples of POTBs, even when their definitions do
not demand it.

    In this paper we present direct and recursive constructions of
POTBs. In the process we have constructed one design  which seems
 to be the first example of an `optimal' two-factor POTB which is not
 a PERGOLA (see Theorem \ref {POTB2}).

\section{Introduction.}

Preece (1966) constructed `BIBDs for two sets of treatments'.
Subsequently several authors constructed similar combinatorial
objects. Among these, the ones relevant to the present paper are
`balanced Graeco-Latin block designs' of Seberry (1979),
`Graeco-Latin designs of type 1' of Street (1981) and `Perfect
Graeco-Latin balanced incomplete block designs (PERGOLAs)' of Rees
and Preece (1999).

 Morgan and Uddin (1996) considered main effects plans (MEPs) on a nested
 row-column set up and proved the optimality of `mutually orthogonal
 balanced nested row-column designs'. They also discussed the constructional
 aspects of such designs. Unfortunately, the relevance of these results in the context of
blocked MEPs was overlooked by later authors studying blocked MEPs
like Mukerjee, Dey and Chatterjee (2001) and Bagchi (2010). Optimal
Main effects plans for three or more factors on non- orthogonal
blocks of small size were obtained in Mukerjee, Dey and Chatterjee
(2001). Bagchi (2010) defined and studied main effect plans
orthogonal `through the block factor' (POTB).

In the present paper we first note the relation between POTBs with
optimality property (termed balanced POTB) and the combinatorial
objects considered by earlier authors mentioned above.
 Next we construct a few series of  POTBs. We also present
a recursive construction by which the number of factors is
multiplied, keeping the block size unchanged, thus yielding
multi-factor POTBs from PERGOLAs and other similar two-factor
designs.

     We note that a PERGOLA is a balanced POTB with an additional
condition. It is interesting that all the balanced POTBs available
in the literature (with different names) do satisfy this condition.
(Table 1 of Rees and Preece (1999) shows that such designs are
plentiful). One would, therefore, suspect that this condition is
implicit in the definition. We have, however, found a balanced POTB
which does not satisfy this condition. [See Theorem \ref {POTB2}].

\section{Preliminaries}
 In this section we present the definition of a balanced POTB. We
 also list related combinatorial objects existing in the literature with
 various names.

\begin{defi} \label{BlockDesign} By a {\bf block design} with
v treatments and b blocks of size k each we mean an {\bf incidence
structure} represented by a $v \times b$ matrix $N$ having
constantcolumn sum $k$.

With any such block design D one associates the graph $G(D)$ with
the treatments of D as vertices, two treatments being adjacent in
$G(D)$ if there is a block containing both the treatments. One says
that {\bf D is connected if the graph $G(D)$ is connected in the
usual sense}.

 $ {\cal D} (b, k, v)$ will denote {\bf the class
of all connected block designs with v treatments on b blocks of size
k each}. \end{defi}

\begin{defi} \label{blockedMEP} Consider an $(m+1) \times n$ array $A$
such that the entries of the ith row are elements of a set ${\mathbf
S}_i$ of size $s_i, \:i=0,1, \cdots m$. This is said to be a {\bf
main effect plan (MEP)} for $m+1$ {\bf factors}, say,  $F_0, F_1,
\cdots F_m$ on $n$ runs. The ith row corresponds to the factor $F_i$
and we say that $F_i$ has $s_i$ levels.

For $0 \leq i,j \leq m$, let ${\mathbf M}_{ij}$ be the $s_i \times
s_j$ matrix such that {\bf the rows and columns of} ${\mathbf
M}_{ij}$ {\bf are indexed by} ${\mathbf S}_i$ and ${\mathbf S}_j$
respectively and {\bf the $(p,q)$th entry of ${\mathbf M}_{ij}$ is
the number of columns of $A$ in which  the ith and jth entries are
$p$ and $q$ respectively}, $ p \in S_i,\: q \in S_j$. ${\mathbf
M}_{ij}$ is said to be {\bf the $F_i$ versus $F_j$ incidence
matrix}.

 Now, suppose {\bf $n =bk$, where  $b$ and $k$ are integers}. A {\bf blocked
 MEP (laid out on b blocks of size k each)} is an $(m+1) \times bk$ array
 in which ${\mathbf S}_0$ is the set of integers $\{1,2 \cdots b \}$
and each integer $j \in S_0$ appears exactly k times in the last row
of $A$. {\bf The 0-th row is said to correspond to the ``block
factor"} , which is represented by $B$ (and not $F_0$). In this
case, {\bf the incidence matrix of the ith row versus the last row
is denoted by} ${\mathbf M}_{iB}, 1 \leq i \leq m$. {\bf A blocked
MEP is said to be `connected'} (borrowing a term from the theory of
block designs) if each ${\mathbf M}_{iB}$ is {\bf the incidence
matrix of a connected block design}. It is said to be {\bf symmetric
if} $s_1 = s_2 \cdots = s_m$. A symmetric MEP with $s_i = s $ for
every $i$ is also referred to as {\bf an MEP for an $s^m$
experiment}.\end{defi}

   In the applications, there are n experimental units, which are
   classified into homogeneous classes or blocks. These units are
   used to study the effects of m factors, the ith one having $s_i$
   `levels'. Typically, an experimental unit, say in the jth block,
   receives a `level combination' say  $x = (x_1, \cdots x_m)$, i.e.
   the level $x_i$ of $F_i$ is  applied on that unit,  $i=1,2, \cdots m$.
   This information is stored in the column vector
   $(j,x_1, \cdots x_m)'$. The array $A$ consists of all such
   column vectors.

\begin{defi} \label{ POTB} [Bagchi(2010)] The ith and jth factors
 of a blocked MEP $\rho$ are said to be {\bf
orthogonal through the block factor (OTB)} if

\begin{equation}  \label{orthbl-inc}
{\mathbf M}_{iB} ({\mathbf M}_{jB})' = k {\mathbf M}_{ij}.
\end{equation}

  If every pair of factors of a plan $\rho$ is orthogonal to each other
  through the block factor, then $\rho$ is said to be a plan orthogonal
through the block factor ({\bf POTB}).
\end{defi}

{\bf Remark 2.1:} What is the utility of condition (\ref
{orthbl-inc}) ? This condition guaranties that  for the inference on
a factor $F_i$ of a POTB one has to look at only its incidence with
the block factor (i.e. $M_{iB}$) and forget all other treatment
factors. Thus, the performance of a POTB $\rho$ regarding the
inference on the ith treatment factor depends only on the incidence
matrix $M_{iB}$.We present a more precise statement in the next
theorem. We omit the proof which can be obtained by going along the
same lines as in the proofs of Lemma 1 and Theorem 1 of Mukerjee,
Dey and Chatterjee (2001). [See Shah and Sinha (1989) for
definitions, results and other details about optimality].

\begin{theo} \label {optimality} Suppose a connected POTB $\rho^*$
satisfies the following condition. For  some non-increasing
optimality criterion $\phi$, $M_{iB}$ is the incidence matrix of a
block design $d$ which is $\phi$ -optimal in the class of all
connected block designs with $s_i$ treatments and b blocks of size k
each. Then, $\rho^*$ is $\phi$-optimal in the class of all connected
m-factor MEP in the same set-up for the inference on the ith
factor..
\end{theo}

In particular, using the well-known optimality property of a BIBD we
get the following result.

\begin{cor}
\label{u-o OMEP} Suppose $\rho^*$ is a connected POTB. Suppose
further that $M_{iB}$ is the incidence matrix of a BIBD. Then,
$\rho^*$ is universally optimal in the class of all m-factor
connected MEP in the same set up, for the inference on the ith
factor. \end{cor}

In view of the above result, we introduce the following term.

\begin{defi} \label{balanced}
A connected POTB is said to be {\bf balanced} if {\bf each of its
factors form a BIBD with the block factor}, that is  $M_{iB}$ is the
incidence matrix of a BIBD for each factor $F_i$.
\end{defi}

We now present a small example of a balanced POTB. This has two
factors, each with four levels 0,1,2,3 on six blocks of size two
each.

\vspace{.5em}

 $ \begin{array}{ccccccccccccccccccc}
\mbox {Blocks}& |& B_1 & & |& B_2 & & |& B_3 & & | & B_4 & & |& B_5 &&|& B_6 &\\
 \hline
          F_1  &|&0 & 2  &|&1  & 3 &|& 0 & 3 &|& 1 & 2      &|&0  & 1 &|&3 &2\\
          F_2  &|&1 & 3  &|&0 &  2 &|& 2 & 1 &|& 3 & 0      &|&3  & 2 &|&0 &1\\
\hline\\\end{array} .$

\vspace{.5em}

 Next we list a few combinatorial designs and note their
 relation with balanced POTBs.

(a) {\bf Balanced  Graco-Latin block design} defined and constructed
in Seberry (1979) are balanced POTBs with two factors.

(b) {\bf  Graco-Latin block design of type 1} of  Street (1981) are
also two-factor  balanced POTBs having

 $$ {\mathbf M}_{12} = J.$$

 (c)  {\bf Perfect Graeco-Latin balanced incomplete block designs
 (PERGOLAs)} defined and discussed extensively in Rees and Preece (1999)
are  two-factor  balanced POTBs having

  (i) $s_1 = s_2 = s$, say and

\begin{equation} \label{extracond} {\mathbf M}_{12} {\mathbf M}'_{12} =
{\mathbf M}'_{12} {\mathbf M}_{12} = f I_s + g J_s, \; \mbox{ f, g
are integers}.
\end{equation}

 (d) {\bf Mutually orthogonal BIBDs} defined and constructed by Morgan and Uddin
 (1996) are multi-factor  balanced POTBs.

 Here $I_n$ is the identity matrix and $J_n$ is the
 all-one matrix of order $n$.

\vspace{.5em}

{\bf Remark 2.2:} The definition of neither  balanced  Graco-Latin
block designs nor of  mutually orthogonal BIBDs include condition
(\ref {extracond}). However, it is interesting to note that all
these designs constructed so far do satisfy this condition. [See
theorem \ref {AllPergola}].

\section{ Constructions for symmetric POTBs}
Now we present a few constructions of POTB's.
 Each of these constructions is in terms of some group $G$ of
order $g$ (which is the additive group of the field $V$ in Theorem
\ref {POTB2new}).

 {\bf A block will consist of k plots or runs represented by
columns. By adding an element $u \in G$ to a block we mean adding
$u$ to the level of every factor in every run of the block.  By
developing an initial block we mean generating g blocks by adding
distinct elements of $G$ to the initial block.}

   Let $N$ denote the set of integers modulo $n$ and $ N^+$ denote $N \cup
\{\infty\}$.

\begin{theo} \label{POTB2} Let $n$ be a positive integer $\geq 5$.

(a) Then there exists a  POTB with three factors $F_0,F_1,F_2$ each
having $n+1$ levels on $ b = 6n$ blocks of size two.

(b) In the case $n = 5$, we get a Balanced POTB, which is {\bf not a
PERGOLA}.
\end{theo}

{\bf Proof :} (a) Let $N^+$ be the set of levels for each factor.
The initial blocks $B_{ij},\; i=1,2,\:j = 0,1,2$ are as follows.
Here addition in the suffix of $F$ is modulo 3.

 $$\begin{array}{cccccccc}
         && \mbox {Block} & B_{1j} & |&  \mbox {Block} & B_{2j}& | \\
     \hline
F_{0+j} & | & \infty        & 0      & |& \infty        & 0 & |\\
F_{1+j} & | &   0           & 1      & |& 0             & 2 & | \\
F_{2+j} & | &  -1           &  1     & |& 1             & 2 &
|\\\hline
\end{array}, \; j=0,1,2.$$

 That the design satisfies the required property follows by
 straightforward verification.

(b) Let $n = 5$. One can verify that the incidence matrices satisfy
the following.

 \begin{eqnarray} \label{5incMat}
 M_{ij} &= & \left [ \begin{array} {cccccc}
0 & 2 & 2 & 2 & 2 & 2  \\
2 & 2 & 2 & 1 & 1 & 2 \\
2 & 2 & 2 & 2 & 1 & 1\\
2 & 1 & 2 & 2 & 2 & 1 \\
2 & 1 & 1 &  2 & 2 & 2 \\
2 & 2 & 1 & 1& 2 & 2 \\\end{array}  \right ],\; i,j = 0,1,2 \mbox{
and }\\
 M_{iB} (M_{iB})' &=& 8I_6 + 2 J_6,\: i=1,2,3.
 \end{eqnarray}

We see that each $M_{iB}$ is the incidence matrix of a BIBD.
 Thus, by definition \ref {balanced} it is a balanced POTB. But
neither of $M_{ij}$s satisfy (\ref {extracond}), as is clear from
(\ref {5incMat}). Thus, the two-factor balanced POTB obtained by
ignoring any one of the factors is not a PERGOLA. $\Box$

\begin{theo} \label {POTBk=2,4} Suppose $n$ is an integer $ \geq 5$.
Then there exists

(a) a  POTB with two $n-$level factors on $2n$ blocks and

 (b) a  POTB with four $n-$level factors on $4n$ blocks of size
 $2$ each.

(c)  We get a balanced POTB  in the case $n = 5$ in series (a).
Further, in the case $n = 10$ in  series (b) we get a POTB which is
E-optimal for the inference on each factor.
\end{theo}

{\bf Proof :} The set of levels of each factor is $N$. Let $a,b \in
N$.

  (a) We present initial blocks $B_1$ and $B_2$ below.

\vspace{.5em}

 $\begin{array}{cccccccc}
    & & \mbox {Block} & B_1 & |&  \mbox {Block} & B_2 & | \\
     \hline
F_1 & | &     a         & -a  &| &    b         &  -b & |   \\
F_2 & | &     b         &  -b & |&    -a         & a & |
\\\hline\end{array}.$

\vspace{.5em}

 (b) We present the initial blocks $B_l, \; l = 1, \cdots 4$ below.

 \vspace{.5em}
 $ \begin{array}{ccccccccccccc}
\mbox {Blocks} & |& B_1&  & |& B_2 &  & |& B_3 & & | & B_4 &  \\
 \hline
          F_1 & | &  0  & a& |&a    &-a& |& 0   &b& | &-b & b \\
          F_2 & | & a   &-a& |&0   &-a & |& -b & b& | & 0 & b \\
          F_3 & | & 0   & b& |& b   &-b& |& -a & 0& | & a & -a\\
          F_4 & | & b   &-b& |& 0   &-b& |& a   &-a&| & a & 0 \\ \hline
\end{array} .$

\vspace{.5em}

That these initial blocks generate POTBs can be verified by
straightforward computation.

\vspace{.5em}

(c) For $n=5$, taking $a=1,b=2$ we get a balanced POTB.

 For $n=10$, we take $a=1$ and $b=3$. Then for every $i =1, \cdots 4$,
 $M_{iB}$ is the incidence matrix of a group divisible design with five
 groups, the jth group being the pair of levels $\{j, j+5\} \:
 j=0,\cdots 4$, satisfying $\lambda_0 = 0$ and $\lambda_1 = 1$. This plan
  is, therefore, E-optimal for the inference on all the four factors by
  Takeuchi (1961). $\Box$

\vspace{.5em}

\begin{theo} \label{oddPOTB}(a) There exists a symmetric POTB with four
$n-$level factors on $4n$ blocks of size $2$ each, whenever $ n \geq
9$. We get a balanced POTB in the case $n=9$.

(b) There exists a  symmetric  POTB with four  factors each having
$n+1$ levels on $6n$ blocks of size $2$ each, whenever $ n \geq 7$.
\end{theo}

{\bf Proof :} (a) The set of levels for each factor is $N$. Let
$a,b,c,d \in N$. The initial blocks  $B_l, \; l = 1, \cdots 4$ are
as follows.

\vspace{.5em}

 $ \begin{array}{ccccccccccccc}
\mbox {Blocks} & |& B_1&  & |& B_2 &  & |& B_3 & & | & B_4 &  \\
 \hline
        F_1  & | & a  & -a & |& b  & -b &|& c & -c &|& -d & d \\
        F_2  & | & b  & -b & |& -a & a  &|& -d & d &|& -c & c \\
        F_3  & | & c & -c  & |& d  & -d &|&-a & a &|&  b  & -b \\
        F_4  & | & d & -d  & |& -c & c &|&  b  & -b&|&a  & -a  \\\hline
        \end{array} .$
\vspace{.5em}

 One can easily verify that these initial blocks generate a symmetric
 POTB with the given parameters.
By taking $a=1,b=2,c=3$ and $d=4$ in the case $n = 9$, we get a
balanced POTB.

\vspace{.5em}

(b) The set of levels for each factor is $N^+$. Let $a,b,c \in N$.
The initial blocks  $B_l, \; l = 1, \cdots 6$ are as follows.

\vspace{.5em}

 $ \begin{array}{ccccccccccccccccccc}
\mbox {Blocks} & |& B_1& & |& B_2 & & |& B_3 & & | & B_4 & & |& B_5 &&|& B_6 &\\
 \hline
  F_1  &| &0 &\infty&|&a  & -a &|& b & -b &|& c & -c &|&a  & -a &|&a &-a\\
  F_2  &| &a  & -a&|&0 &\infty &|& c & -c &|& -b & b &|&a  & -a &|&-a &a\\
  F_3  &| & b & -b &|& c & -c  &|&0 &\infty&|&a  & -a&|&-c & c &|&-c & c\\
  F_4  &| &c & -c &|&b & -b    &|&a & -a   &|&0&\infty&|&-c & c &|&c &-c\\
\hline\\\end{array} .$

\vspace{.5em}

That the design satisfies the required property follows by
 straightforward verification. $\Box$

 Next we construct a series of balanced POTBs using finite fields.
 We first  introduce the following notation.

\begin{nota} \label{Unionnota}

(i) Let  $\bigsqcup$ denote an union counting multiplicity.

(ii) For a set $A$ and an integer $n$, let $nA$ denotes a multiset
in which every member of $A$ occurs $n$ times.

(iii) For subsets A and B of a group $(G,+)$,
$$ A - B = \{ a-b : a \in A, b \in B \}.$$
\end{nota}

\begin{nota} \label{GFnota}(i) v denotes an odd prime or a prime
power, written as $v = m f + 1$. $V$ denotes the Galois field of
order v. Further, $V^* = V \setminus \{0 \}$ and $V^+ = V \cup
\{\infty\}$.

(ii) $\alpha$  denotes a primitive element of $V$.

 (iii)  $\beta = \alpha^m$ is a generator of the subgroup $C_0$ of
 order f of $(V^*,.)$.

(iv) $C_0, C_1, \cdots C_{m-1}$ are the cosets of $C_0$ in
$(V^*,.)$.

(v) $(i,j) =$ the number of ordered pairs of integers (s,t) such
that the following equation is satisfied in $V$. [ This notation is
borrowed from the theory of cyclotomy]

 $$ 1 + \alpha^s = \alpha^t ,\; s \equiv i, t \equiv j \pmod m.$$
 \end{nota}

 We need the following well-known result. [See  Hall
(1986), for instance].

\begin{lem} \label{oddEven} Suppose $m=2$.Then
the following hold.

(a) $-1 \in C_0$ (respectively $C_1$) if f is even
  (respectively odd).

 (b) If $f$ is even, then  $\alpha -1 \in C_i \Rightarrow \alpha^{-1} -1 \in
 C_{i+1}$.

(c) If $f$ is odd, then  $\alpha -1 \in C_i \Rightarrow 1 -
\alpha^{-1} \in  C_{i+1}$.

Here $+$ in the suffix is modulo 2. \end{lem}

 We present the following well-known results for ready
reference. [See equations (11.6.30), (11.6.40) and (11.6.43) of Hall
(1986)].

\begin{lem} \label{cyclotomy}
 The  differences between the cosets $C_j$'s of $V^*$ can be
expressed in terms of the cyclotomy numbers (p,q)'s as follows.

$$ C_i - C_j = \left \{ \begin{array}{ll}
\bigcup\limits_{k=0}^{m-1} (k-j,i-j) C_k & \mbox{ if } j \neq i \\
   f \{0 \} \cup \bigcup\limits_{k=0}^{m-1} (k-j,0) C_k & \mbox{ if } j = i
\end{array}
\right . $$

The following cyclotomy numbers are known for the case $m = 2$.
\end{lem}

{\bf Case 1: f odd.} (0,0) = (1,1) = (1,0) = (f-1)/2, (0,1) =
(f+1)/2.

{\bf Case 2: f even.} (0,0) = f/2 -1, (0,1) = (1,0) = (1,1) = f/2.

\vspace{.5em}

{\bf A series of two-factor balanced POTBs :}

\begin{theo} \label{POTB2new} Suppose $v$ is an odd prime or a prime
power. Then there exists a balanced POTB for a $(v+1)^2$ experiment
on $b = 2v$ blocks of size $(v+1)/2$ . \end{theo}

{\bf Proof :} We write $v = 2f+1$. The set of levels of each factor
is $V \cup \{\infty \}$.  The plan is obtained by developing the
following initial blocks $B_0$ and $B_1$ presented below.

 {\bf Case 1 : f is even.}

$$B_0 = \left [\begin{array}{cccccc}
  \infty & 1      & \beta      & \cdots & \beta^{f-1} \\
  0      & \alpha &\alpha\beta& \cdots & \alpha\beta^{f-1} \end{array}
 \right ] \mbox{ and } B_1 =\left [\begin{array}{cccccc}
  0    & 1      & \beta      & \cdots & \beta^{f-1} \\
 \infty &\alpha^{-1} &\alpha^{-1}\beta& \cdots & \alpha^{-1}
\beta^{f-1} \end{array} \right ].$$

{\bf Case 2 : f is odd.}

Block $B_0$ is as in case 1, while  Block $B_1$ is as follows.
$$ B_1 =\left [\begin{array}{cccccc}
 0 &\alpha^{-1} &\alpha^{-1}\beta& \cdots & \alpha^{-1}
\beta^{f-1} \\
 \infty    & 1      & \beta      & \cdots & \beta^{f-1}
 \end{array} \right ].$$

Clearly block size is $f+1 = (v+1)/2$. To show that the plan
satisfies the required property, we have to show that

(a) the plan is a POTB and (b) each factor forms a BIBD with the
block factor.

Condition (b) follows from Lemma \ref {cyclotomy}. So, we prove (a)
.

let us use the following {\bf simplified notation} ${\mathbf M}
=((m_{ij}))$ for ${\mathbf M}_{12}$ and ${\mathbf A}$ for ${\mathbf
M}_{1B} ({\mathbf M}_{2B})'$. We note that $m_{ij}$ is the total
number of plots (runs) in which the level combination $(i,j)$
appears, while $a_{ij}$ is the number of blocks in which $F_1$ is at
level i and $F_2$ at level j, (in same or different plots).

We shall show that
\begin{eqnarray} \label{RelIncMat}
M &= & J - I \mbox{ and }\\
A &= & (f+1)( J - I) \end{eqnarray}

We begin with $M$. It is clear that
 $ m_{ii} = 0, i \in V^+ $ and
$ m_{\infty,i} =  m_{i,\infty} =1, i \in V$.

\vspace{.5em}

We, therefore, assume $i \neq j,\; i,j \in V$. Let $u = j-i$. Then,
$m_{ij}$ is the number of times $u$ appears in the multiset

$$ \left \{\begin{array} {ll}
(\alpha -1) C_0 \bigsqcup (\alpha^{-1} -1) C_0 & \mbox{ if f is
even}\\
(\alpha -1) C_0 \bigsqcup (1 - \alpha^{-1}) C_0 & \mbox{ if f is
odd} \end{array} \right . $$

The relations above imply (\ref {RelIncMat}) in view of Lemma \ref
{oddEven}.

Now we look at $A$. Clearly, $A_{ii} = 0, \; i \in  V^+$. Further,
for every $i \in V$, $m_{\infty,i}$  is the replication number of
$i$ in the block design generated by the initial block $\{0\} \cup
C_1$. Similarly, $m_{i,\infty}$ is the replication number of $i$ in
the block design generated by the initial block $\{0\} \cup C_0$ if
$f$ is odd and $\{0\} \cup C_1$ otherwise. Thus,

$$m_{\infty,i} = m_{i, \infty} = f+1, \; i \in V.$$

 Next we consider $i \neq j, \; i,j \in V.$ Let $u = j-i$. Then, $a_{ij}$
is the number of times $u$ appears in the multiset

$$ \tilde{S} =  \left \{\begin{array} {ll}
((\{0\} \cup C_1) - C_0) \bigsqcup (C_1 - (\{0\} \cup C_0)) & \mbox{
if f is even}\\
((\{0\} \cup C_1) - C_0) \bigsqcup (C_0 - (\{0\} \cup C_1)) & \mbox{
if f is odd}.\end{array} \right . $$

These, together with  Lemma \ref {cyclotomy} and (a) of Lemma \ref
{oddEven} imply the equation next to  (\ref {RelIncMat}). $\Box$

\vspace{.5em}

 Now we present the series of balanced POTBs available in the
 literature, together with two newly constructed balanced POTBs.

\begin{nota}  v denotes an odd prime or a prime power of the form
$v = mf + 1$. \end{nota}

\begin{center}
{\bf {\large Table 3.1 : Balanced POTBs}}

\vspace{.5em}

\begin{tabular}{ccccccc} \hline
 No. & The expt. & m & \# of Blocks & Block size & Inc. matrix ($M_{ij})$ &
 Reference\\
 \hline\\
 1. &$(v+1)\times v$ & 2& $2v$& $f =(v-1)/2$ &  $J$& Seberry (1979)\\
 2. & $(v+1)^2$ & 2& $2v$& $f +1 =(v+1)/2$ & $ J$&  Street (1981) \\
 3(a). & $v^m$  & $m$&$mv$& $f$ & $ (J - I)$&Morgan
  and Uddin(1996)\\
 3(b). & $v^t$ & $tg$&$mv$& $hf,\: h \leq g$& $ h(J - I)$&Morgan
  and Uddin(1996)\\
 4. & $v^f$    & m &$mv$& $1 + hf,\: h \leq m$&$(m-h)I + hJ$&Morgan
 and Uddin(1996)\\
 5. & $(v+1)^2$ &2   & $2v$&$f+1$   &$ J - I$& Theorem \ref
 {POTB2new}\\
 6. &$(5+1)^3$&- &30 &2  & given in(\ref {5incMat}) & Theorem \ref {POTB2}\\
 7. &$9^4$& - & 36 & 2 & $ J - I$ & Theorem \ref {oddPOTB}\\
\hline\end{tabular}  \end{center}

    In 3(b) above $t$ is a factor of $m, g = m/t$ and $h$ is an
    integer $ \leq g$.

\vspace{1em}

 Using the information in Table 3.1, one may verify the following
 result.

\begin{theo} \label{AllPergola} Every two-factor balanced POTB obtained from an existing
multi-factor balanced POTB, except the one with $n=5$, constructed
in Theorem \ref {POTB2}, is a PERGOLA.\end{theo}

One may also look at  Table 1 of Rees and Preece (1999) for many
more examples of PERGOLAS.

\vspace{.5em}

 {\bf A recursive construction }

\begin{nota} An orthogonal array with $m$ rows, $n$ columns, k symbols and
strength 2, will be denoted by {\it $OA(n,m,k,2)$}. \end{nota}

\begin{theo} Suppose there exists a balanced POTB with $f$ factors on $b$
blocks  of size $k$ each, with $f \leq k$. If further  an $OA (k^2,
m+1,k,2)$ exists, then a  balanced POTB with $mf$ factors on $bk$
blocks of size $k$ each also exists.
\end{theo}

The proof of this theorem is based on the following lemma.

\begin{lem}\label{useOrth}  Consider a set of $k$ runs of a plan for an experiment
with $f (\leq k)$ factors, such that no level of any factor is
repeated. If  an $OA (k^2, m+1,k,2)$ exists, then there exists an
MEP with $f$ classes of m $k-$level factors on $k$ blocks  of size
$k$ each with the following property. Every factor is orthogonal
(w.r.t. the block factor) to  every factor of a different class.
\end{lem}

 {\bf Proof :} Let $D$ denote the given set of runs. Let $F =\{ P,Q,
 \cdots \}$ denote the set of f factors of $D$. For each $P \in F$,
 let $K_P$ denote the set of levels of P appearing in $D$.  Let $K = \{1, \cdots
 k \} $ denote the set of symbols of the given OA $(\tilde{O})$, say.
 For every $P \in F$, let $L_P$ denote the following one-one function
 from $K$ to $K_P$.  [By assumption, size of $K_P$ is $k$ for each $P$].

 \begin{equation}  \label{Level-corr} L_P(i) = j,\; i \in K, j \in K_P,
 \mbox{ if P has level j in the ith run of D}.\end{equation}

 Let us arrange the columns of the given OA $(\tilde{O})$ as

 $$ (\tilde{O}) = \left [ \begin{array}{ccc}
\tilde{A_1} & \cdots \tilde{A_{k}}
\end{array} \right ], $$

such that the 1st row of $\tilde{A_i}$ consists of the symbol $i$
repeated k times,  $1 \leq i \leq k$. Let $A_i$ denote the $m \times
k$ array obtained from $\tilde{A_i}$ by deleting the 1st row. Thus,
every member of $K$ appear exactly once in every row of each $A_i, i
=1, \cdots k$.

We now construct $D^*$, the reqd MEP. For each factor $P$ of D,
there will be m factors $P_1, \cdots P_m$ in $D^*$, each of which
will have $K_P$ as the set of levels.

For $i \in K$ , let us fix $A_i$ and a factor, say $P$ of $D$. If
the $j$ th column of $A_i$ is $(s_1, \cdots s_m), s_u \in K$, then
in the j-th plot of the $i$-th block of $D^*$, the factor $P_t$ will
have level $L_P(s_t), t=1,2, \cdots m$, where $L_P$ is as in (\ref
{Level-corr}). Doing the same for all the factors and varying $j$
over $ \{1,2, \cdots k\}$ we get a block of $D^*$. Finally varying
$i$ over $K$ we we generate the k blocks of the reqd MEP.

We now show that the MEP $D^*$ satisfies the required property. We
fix two factors, say $P_i$ and $Q_j, \; i \neq j$ and an ordered
pair of levels, say $(u,v), u \in K_P,\: v \in K_Q$. From the
construction the following is clear. In every block there is a plot
in which $P_i$ is at level u  and a plot where $Q_j$ is at level v.
Moreover, there is exactly one block in which these factors are set
at these levels in the same plot. Thus, the factors $P_i$ and $Q_j$
are mutually orthogonal through the block factor. We see that if
$P=Q$, then also the argument above holds. Thus $P_i$ is orthogonal
to $P_j, j \neq i$. However, $P_i$ may not be orthogonal to $Q_i$.
We, therefore form the classes as $C_i = \{ P_i, Q_i, \cdots  \}, \;
P, Q \in F, \: i = 1,2, \cdots m$. Now the factors satisfy the
orthogonality condition of the hypothesis.$\Box$

{\bf Proof of the theorem :} Let $D$ denote the given POTB. Let
$A_i, i \in K$ be as in Lemma. For every block of $D$ we construct a
an MEP following the method described in the proof of the lemma
above. Let the resultant MEP be named $D^*$. By Lemma \ref
{useOrth}, every pair of factors belonging to different classes are
orthogonal w.r.t. the block factor. Further, since the pair of
factors $P,Q$ are mutually orthogonal w.r.t. the block factor in
$D$, it follows that for every $i \in K$, the factors $P_i$ and
$Q_i$ are also  mutually orthogonal w.r.t. the block factor in
$D^*$. $\Box$

{\bf Remark 3.1:} If we look at the restriction of $D^*$ to one
factor, say $P_i$, we see that it is nothing but k times repetition
of each block of the block design obtained from the restriction of
$D$ to the factor P.

\vspace{.5em}

\section{References}
\begin{enumerate}
\item  Bagchi, S. (2010). Main effect plans orthogonal through the
block factor. Technometrics, vol. 52, p. 243-249.

\item Hall, M. (1986). Combinatorial Theory. Wiley-interscince, New
York.

\item Morgan, J.P. and Uddin, N. (1996). Optimal blocked main effect
plans with nested rows and columns and related designs. Ann. Stat.
vol. 24, p. 1185-1208.

\item Mukerjee, R., Dey, A. and Chatterjee, K. (2001). Optimal
main effect plans with non-orthogonal blocking. Biometrika, 89, p.
225-229.

\item  Preece, D.A. (1966). Some balanced incomplete block designs
for two sets of treatment. Biometrika 53, p. 497-506.

\item Rees, D.H. and Preece, D.A. (1999). Perfect Graeco-Latin
balanced incomplete block designs. Disc. Math. vol.197/198, p.
691-712.

\item Seberry, Jennifer, (1979). A note on orthogonal Graeco-Latin
designs. Ars. Combin. vol. 8, p. 85-94.

\item Shah, K.R. and Sinha, B.K. (1989). Theory of optimal
designs, Lecture notes in Stat., vol. 54, Springer-Verlag, Berlin.

\item Street, D.J. (1981). Graeco-Latin and nested row and column
designs. In Com. Math. VIII, Proc. 8th Austr. Conf. Comb. Math.,
Lecture notes in Math., vol. 884, Springer, Berlin. p. 304-313.

\item Takeuchi, K. (1961). On the optimality of certain type of PBIB
designs. Rep. Stat. Appl. Un. Jpn. Sci. Eng. vol. 8. p. 140-145.
\end{enumerate}

\end{document}